\documentclass[a4paper,12 pt]{article}
\usepackage[utf8]{inputenc}
\usepackage[margin=1 in]{geometry}

\usepackage{siunitx}

\usepackage{amsthm}
\usepackage{amsmath}
\usepackage{amsfonts}
\usepackage{dutchcal}
\usepackage{amssymb}
\usepackage{epsfig,epstopdf,titling}
\usepackage{caption}
\usepackage{subcaption}
\usepackage{multirow}
\usepackage{url}
\usepackage{graphicx}
\usepackage{pgfplots}
\usepackage{array,tabularx}
\usepackage{multicol}  
\usepackage{algorithmic} 
\usepackage[boxruled,longend,linesnumbered,ruled]{algorithm2e}
\usepackage{longtable}
\usepackage{amssymb}
\usepackage{multirow}
\usepackage{cite}
\usepackage{hyperref}
\usepackage{tikz}
\usetikzlibrary{positioning}
\usetikzlibrary{calc}
\usepackage{xcolor}
\usetikzlibrary{shapes.geometric, arrows}
\usepackage[export]{adjustbox}
\usepackage{float}
\usepackage{bm}

\tikzset{mark options={solid,mark size=3,line width=1pt,mark repeat=15}}
\pgfplotsset{every axis plot/.append style={line width=1.5pt}}

\title{\textbf{A survey of the noise-correcting tools for Dynamic Mode Decomposition}}

\author{Moajjem H. Chowdhury\thanks{Department of Electrical and Computer Engineering, North South University} \and Md. Nazmul Islam Shuzan\footnotemark[1]
\and Mohammad N. Murshed\thanks{Department of Mathematics and Physics, North South University} \and Sanwar Alam\footnotemark[1] \and M. Monir Uddin\footnotemark[2] \and Zarin Subah\thanks{Institute of Water and Flood Management, Bangladesh University  of Engineering and Technology}}
\begin{document}

\maketitle
\begin{abstract}
Dynamic Mode Decomposition (DMD) is a data-driven modeling tool that generates a model from spatio-temporal data. The data needs to be as clean as possible for DMD to come up with a faithful model. We review a few data-filtering methods to be integrated with DMD and test them on datasets of varying complexity. The impact of SNR on these methods and the error variation in the DMD model due to each method are observed and discussed.
    
\end{abstract}

\section{Introduction}
There is a lot of data relevant to what is going on around us, be it from the ocean-atmosphere interaction or the turbulent blood flow in the veins. Problems vary from being simple (linear) to complex (highly non-linear). But, one common aspect in the data is the noise. Noise can be of two types: white and colored. The difference between white noise and colored noise is that the power spectral density for the latter varies with frequency. Noise in the data is an issue when making a reduced order model based on that dataset. Reduced order models (ROM) are approximate models of a given full order model (FOM). ROMS appear to be very effective in terms of time and storage. A few notable reduced order models are Proper Orthogonal Decomposition \cite{chatterjee2000introduction}, Balanced Proper Orthogonal Decomposition \cite{willcox2002balanced}, Eigensystem Realization Algorithm \cite{juang1985eigensystem,murshed2020towards}, and Dynamic Mode Decomposition (DMD) \cite{jtu,schmid2010dynamic,9038561,murshed2020projection}.\\ \\
DMD is a data-driven method of creating a model solely from time-series data whose performance usually declines when dealing with corrupted data. There are two types of noise: process noise and sensor noise. The former relates to the system dynamics and latter to the measurements. There are several variants of DMD, namely Total least square DMD \cite{hemati2017}, Forward-Backward DMD \cite{dawson2016characterizing}, and Optimized DMD \cite{askham2018variable}, which are developed to tackle the effect of stochastic sensor disturbance on the dataset. \\ \\
The rest of the paper is organized as follows: Section \ref{BACK} is a brief discussion of the basics of the DMD and some of the remarkable data-filtering recipes. The results are shown in Section \ref{NR} and a summary provided in Section \ref{CFW}. 

\section{Background}
\label{BACK}
This section comprises a review of the algorithm behind Dynamic Mode Decomposition and a precise discussion about the data-filtering methods that would be considered as a pre-processing step in DMD.
\subsection{Dynamic Mode Decomposition}
Given time-series data, Dynamic Mode Decomposition approximates the given data and predicts approximate states in future. The idea of extraction of dynamic information from snapshots is closely related to the Arnoldi algorithm \cite{arnoldi1951principle}. If we have available \(P\) number of snapshots, then, we consider those snapshots in a matrix 
\begin{equation}
    \textbf{X}_{1}^{P} = [\textbf{x}_{1}\ \textbf{x}_{2}\ ...\ \textbf{x}_{P}],
\end{equation}
where each snapshot contains \(Q\) pixels. It is assumed that each snapshot is separated by some uniform time interval, \(\Delta t\). DMD is thought to be the approximation of \emph{Koopman operator} (a linear, infinite-dimensional operator that represents nonlinear, infinite-dimensional dynamics). DMD algorithm considers a best fit linear operator \textbf{A} to approximate the dynamics,
\begin{equation}
    \textbf{x}_{k+1} \approx \textbf{A}\textbf{x}_k,
    \label{equation_3}
\end{equation}
where \(\textbf{x}_k\) is the snapshot at time \(t_k\) and \(\textbf{x}_{k+1}\) the snapshot one time step ahead in future. The data in the matrix can also be written, in terms of the operator, as
\begin{equation}
    \textbf{X}_{2}^{P} \approx \textbf{A}\textbf{X}_{1}^{P-1}.
\end{equation}
The eigendecomposition of \textbf{A} facilitates analysis of the data despite the large size of the operator matrix.\\ \\
\textbf{Algorithm of DMD}\\ \\
(1) DMD takes the singular value decomposition (SVD) of \textbf{X} \cite{trefethen1997numerical}: 
\begin{equation}
    \textbf{X}_{1}^{P-1} = \textbf{U}\Sigma \textbf{V}^\ast
    \Rightarrow \textbf{X}_{2}^{P} = \textbf{A} \textbf{U}\Sigma \textbf{V}^\ast,
\end{equation}
here matrix \(\textbf{U}\) is \(q \times r\), \(\Sigma\) is \(r \times r\) and \(\textbf{V}\) is \(m \times r\). \(\textbf{V}^\ast\) denotes the conjugate transpose of \textbf{V} and \(r\) is the rank for truncation after SVD. Note that the columns in \(\textbf{U}\) represent the POD modes. \\ \\
(2) \(\tilde{\textbf{A}}\) is computed as the \(K \times K\) projected version of matrix \(\textbf{A}\) in POD modes:
\begin{equation}
    \textbf{A} = \textbf{X}_{2}^{P} \textbf{V}\Sigma^{-1} \textbf{U}^\ast
    \Rightarrow \tilde{\textbf{A}} = \textbf{U}^\ast \textbf{A} \textbf{U} = \textbf{U}^\ast \textbf{X}_{2}^{P} \textbf{V}\Sigma^{-1}.
\end{equation}
(3) The eigendecomposition of \(\tilde{\textbf{A}}\) leads to,
\begin{equation}
    \textbf{A}\textbf{W} = \textbf{W}\Lambda,
\end{equation}
where the columns of \(\textbf{W}\) represents the eigenvectors and \(\Lambda\) is a diagonal matrix that contains the corresponding eigenvalues.\\ \\
(4) DMD model reads,
\begin{equation}
    \textbf{X}_{DMD} (t) = \bm{\Phi} \text{exp}(\bm{\Omega} t))\textbf{b},
    \label{equation_9}
\end{equation}
where \(\bm{\Phi} = \textbf{X}_{2}^{P} \textbf{V}\Sigma^{-1} \textbf{W}\) and \(\textbf{b} = \bm{\phi}^{\dagger}\textbf{x}_1\). \(\bm{\Phi}\) is a matrix containing the eigenvectors, and \textbf{b} the vector where each entry is the initial amplitude of each mode. Note that the dagger symbol used in the definition of \textbf{b} denotes the Moore-Penrose pseudoinverse.
\subsection{Noise in data}
We have just seen how DMD begins working with the data matrix, \textbf{X}, to form a model for the dynamics of the data. In practice, the dataset would pick up noise from many sources. In case of numerical data, improper step size in time and space may introduce error in the simulation. For experimental data, uncertainties arise due to poor quality of the instrument and the way the technician is using it. We define a corrupted snapshot as,
$$ \textbf{x}_{c}(t) = \textbf{x}(t) + \bm{\eta},$$
where \(\bm{\eta}\) is a vector containing white Gaussian noise. The first subset of data can be written as,
$$ \textbf{X}_{c} =\textbf{X} + \textbf{C}_{X},$$
where \(\textbf{X}_{c}\) is the corrupted data, \(\textbf{X}\) the noise-free data, and \(\textbf{C}_{X}\) the noise. In the same manner, the second subset would appear as,
$$ \textbf{Y}_{c} =\textbf{Y} + \textbf{C}_{Y}.$$
Taking into account the error in the subsets, the linear mapping that takes place in DMD then becomes,
$$ \textbf{A} = (\textbf{Y}_{c} + \textbf{C}_{Y}) (\textbf{X}_{c} +  \textbf{C}_{X})^{\dagger}, $$
$$ \textbf{A} = (\textbf{Y}_{c} + \textbf{C}_{Y}) (\textbf{X}_{c} +  \textbf{C}_{X})^{*} \{(\textbf{X}_{c} +  \textbf{C}_{X})(\textbf{X}_{c} +  \textbf{C}_{X})^{*}\}^{\dagger}. $$
\subsection{Data-filtering methods}
Many algorithms have been proposed to remove noise from corrupted data. Possible ways are Principal Component Analysis (PCA) and Robust Principal Component Analysis (RPCA)\cite{wright2009}.  RPCA handles noise of large magnitude much better than PCA. RPCA is based on the idea that the data, \textbf{X}, can be decomposed into two matrices: \(\textbf{X} = \textbf{L} \ + \ \textbf{S}\), where \textbf{L} is a low-rank matrix and \textbf{S} a sparse matrix. The idea of calculating \textbf{L} and \textbf{S}, hence, becomes a minimization problem. Previously, this has been done by Iterative Threshold\cite{wright2009}, Accelerated Proximal Gradient \cite{lin2009}, Dual Approach\cite{lin2009}, Singular Value Thresholding \cite{cai2010}, Alternating Direction Method \cite{yuan2009}, Exact Augmented Lagrange Multiplier \cite{lin2010}, and Inexact Augmented Lagrange Multiplier \cite{lin2010}. Variants of DMD can also be used to clean the data. Forward-Backward DMD and Total Least Squares DMD (tlsDMD) \cite{hemati2017} are a few notable ones that complete such task. In what follows, we will review the noise reduction by
\begin{itemize}
\item RPCA with ADM,
\item RPCA with inexact ALM, 
\item Total least-squares DMD.
\end{itemize}
\textbf{RPCA with ADM}\\ \\
RPCA with ADM algorithm uses two helper functions: Shrink and Singular Value Threshold (SVT).\\ \\
\textbf{Shrink}: it takes two parameters as input: \textbf{X} and \(\tau\), where \textbf{X} can be any matrix and \(\tau\) is the threshold. It first checks each element in \textbf{X} and compares its absolute value with \(\tau\). If the absolute value of the element is greater \(\tau\), then the element is replaced with the difference between its absolute value and \(\tau\), otherwise, it is replaced with zero. The resulting matrix is then multiplied by the sign of the original matrix \textbf{X}. The final product is then returned as an output.\\ \\
\textbf{SVT}: this function takes as inputs \textbf{X} and \(\tau\) (a threshold for the consideration of the singular values). The matrix \textbf{X} is first decomposed into its component using SVD. We will then use shrink operator on the singular value matrix. This will cause many small values in the singular value matrix to become zero. The product of the \textbf{U}, \textbf{V}, and the modified \textbf{S} matrix is returned. The implication of this is that the \textbf{X} matrix will become more low-ranked.\\ \\
\noindent \textbf{RPCA with inexact ALM}\\ \\
The augmented Lagrange multiplier (ALM) is a class of algorithm for solving constrained optimization problems. ALM is used for problems where we have to minimize $f(X)$ such that $h(X)=0$. The ALM can be then defined as:
\begin{center}
$L(\textbf{X},\textbf{Y},\mu)=f(\textbf{X})+<\textbf{Y},h(\textbf{X})>+\frac{\mu}{2}{\lvert\lvert h(\textbf{X}) \rvert\rvert}^2_F. $
\end{center}
This algorithm is very fast in solving RPCA compared to other methods such as the Accelerated Proximal Gradient approach.
For the RPCA problem,
\begin{center}
    $ \bm{X} = (\textbf{A,E})$,\\
    $f(\textbf{X}) = \textbf{A}_* + \lambda{\lvert\lvert \textbf{E} \rvert\rvert}_1$,\\
    $h(\textbf{X}) = \textbf{D}-\textbf{A}-\textbf{E}$
\end{center}
such that the lagrangian becomes,
\begin{center}
    $L(\textbf{A,E,Y},\mu) = \textbf{A}_* + \lambda{\lvert\lvert \textbf{E} \rvert\rvert}_1+<\textbf{Y},\textbf{D}-\textbf{A}-\textbf{E}>+\frac{\mu}{2}{\lvert\lvert \textbf{D}-\textbf{A}-\textbf{E}} {\rvert\rvert}^2_F ,     $
\end{center}
where \(\lambda = \frac{\lambda_0}{\sqrt{max(n,m)}}\) with \(\lambda_0\) optimally being 1 and \(\mu\) is a hyperparameter. The inexact ALM can then be calculated as prescribed in \cite{lin2010}.\\ \\
\textbf{Total least-squares DMD}\\ \\
Data is collected and separated into \(\textbf{X}_1\) and \(\textbf{X}_2\). We define \textbf{Z} as,
\begin{equation}
 \textbf{Z} = \textbf{X}_1^{'}\textbf{X}_{1} + \textbf{X}_{2}^{'}\textbf{X}_{2}.
\end{equation}
Then, SVD is performed on \textbf{Z},
\begin{equation}
\textbf{Z}=\textbf{U}\textbf{S}\textbf{V}^*,
\end{equation}
to identify the number of significant modes, \(r\).
$\textbf{V}$ is truncated to $\textbf{V}_n$, 
$$ [\textbf{V}_{n}, \textbf{D}] = eig(\textbf{Z}) $$
$$ \textbf{d} = diag(\textbf{D}) $$ 
$$ [ \sim , idx] = sort(\textbf{d}, 'descend') $$
$$ \textbf{V}_{n} = \textbf{V}_{n}(:,idx(1:r)),$$
to get projected \(\textbf{X}_1\) and \(\textbf{X}_2\):
\begin{center}
$\tilde{\textbf{X}_1}=\textbf{X}_{1} \textbf{V}_{n} \textbf{V}_n^*,$\\
$\tilde{\textbf{X}_2}=\textbf{X}_{2}  \textbf{V}_{n} \textbf{V}_{n}^*.$
\end{center}
Note that \textbf{eig} means the eigendecomposition of a matrix, \textbf{diag} picks out the elements of a matrix from the main diagonal, and \textbf{sort} arranges the entries in the given vector in either ascending or descending order.\\ \\
SVD on $\tilde{\textbf{X}_1}$ yields,
\begin{center}
    $\tilde{\textbf{X}_{1}}=\tilde{\textbf{U}}\tilde{\textbf{S}}\tilde{\textbf{V}}.$
\end{center}
$\tilde{\textbf{A}}$ is then computed as,
\begin{center}
    $\tilde{\textbf{A}} = \tilde{\textbf{U}}\tilde{\textbf{X}_2}\tilde{\textbf{V}}\tilde{\textbf{S}}^{-1}.$
\end{center}

\section{Numerical Results}
\label{NR}
The discussed filtering techniques will be implemented on different dataset. Gaussian white noise is added to the raw data to mimic the corrupted data. The experiments are run on MATLAB 2019b on a laptop with configuration of Intel Core i7-7700HQ CPU @ 2.80 GHz 2.81 GHz with 16 GB RAM.\\ \\
The performance of the model is to be evaluated as per two measures: Root Mean Squared Error (RMSE) and Correlation Co-efficient (CC). These two metrics are defined below, where, \(\textbf{X}_{DMD}\) is the prediction, \textbf{X} the ground truth, and \(n_s\) the number of samples or recordings.
\begin{itemize}
\item RMSE measures the standard deviation of the prediction error (residuals) as,
\begin{equation}\label{rmse}
    RMSE = \sqrt{\frac{\sum {\lvert \textbf{X}_{DMD} - \textbf{X} \rvert}^2}{n_s}}
\end{equation}
The smaller the spread, the better the model.
\item CC indicates how closely two variables (prediction and ground truth) are related and is computed as,
\begin{equation}\label{cc}
    CC = \sqrt{1 - \frac{\frac{\sum \lvert \textbf{X}_{DMD} - mean(\textbf{X}_{DMD}) \rvert}{n_s}}{\frac{\sum \lvert \textbf{X} - mean(\textbf{X}) \rvert}{n_s}}}.
\end{equation}
CC allows us to understand if the general structure of the signal/data is retained.
\end{itemize}
The error in the DMD approximation over time will be computed as,
$$ \epsilon = \frac{|\textbf{x}_{DMD} - \textbf{x}|_{2}}{|\textbf{x}|_{2}}, $$
where \textbf{x} is the approximation at a certain time instant.

\subsection{Non-linear Schrodinger Equation}
\begin{figure}[H]
    \centering
    \begin{subfigure}{0.47\textwidth}
    \begin{adjustbox}{max width=1\textwidth,center}
    \includegraphics{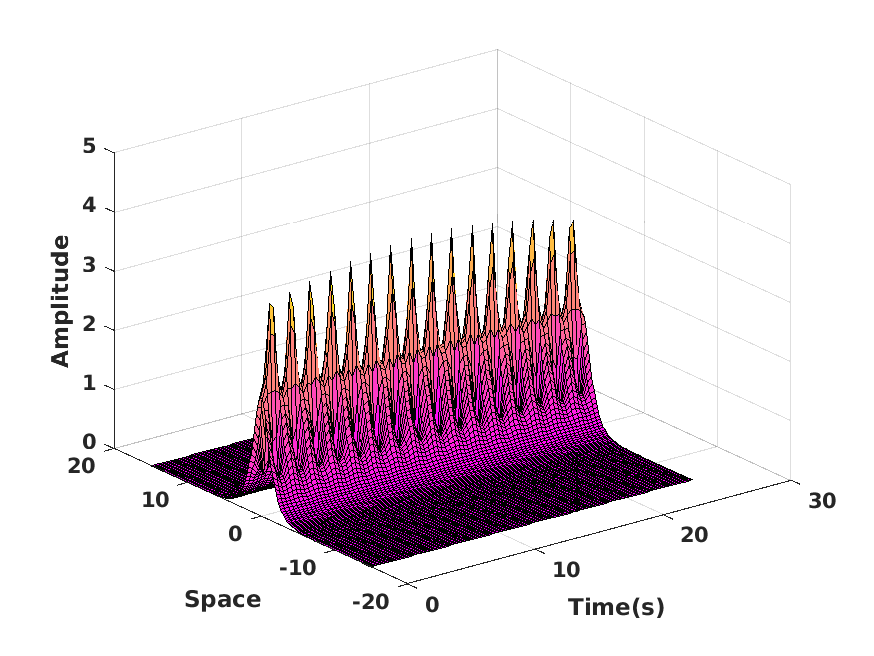}
    \end{adjustbox}
    \caption{Original Data}
    \label{nlse_odat}
    \end{subfigure}
    \hspace{0.04\textwidth}
    \begin{subfigure}{0.47\textwidth}
    \begin{adjustbox}{max width=1\textwidth,center}
    \includegraphics{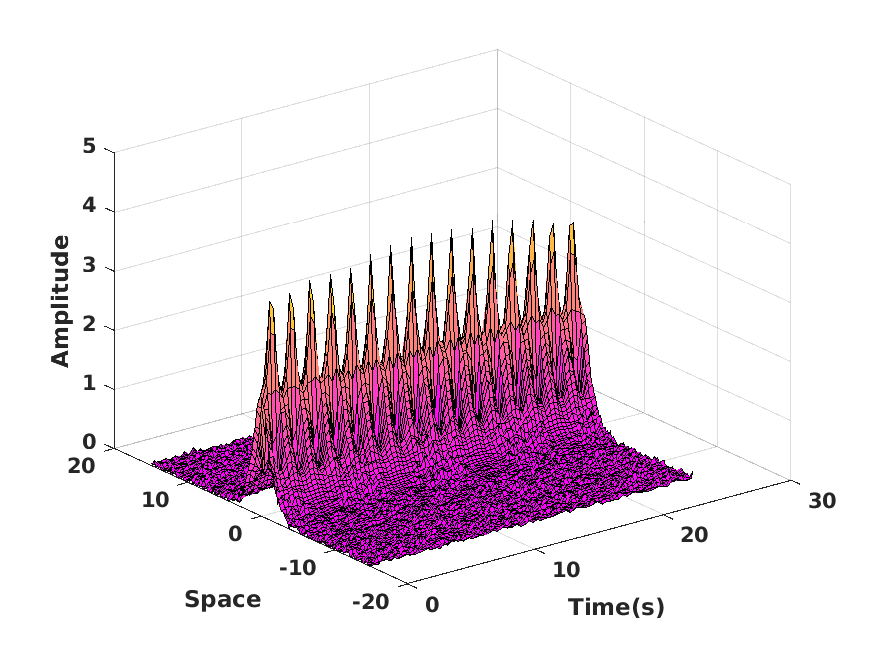}
    \end{adjustbox}
    \caption{Noisy Data}
    \label{nlse_ndat}
    \end{subfigure}
    \begin{subfigure}{0.47\textwidth}
    \begin{adjustbox}{max width=1\textwidth,center}
    \includegraphics{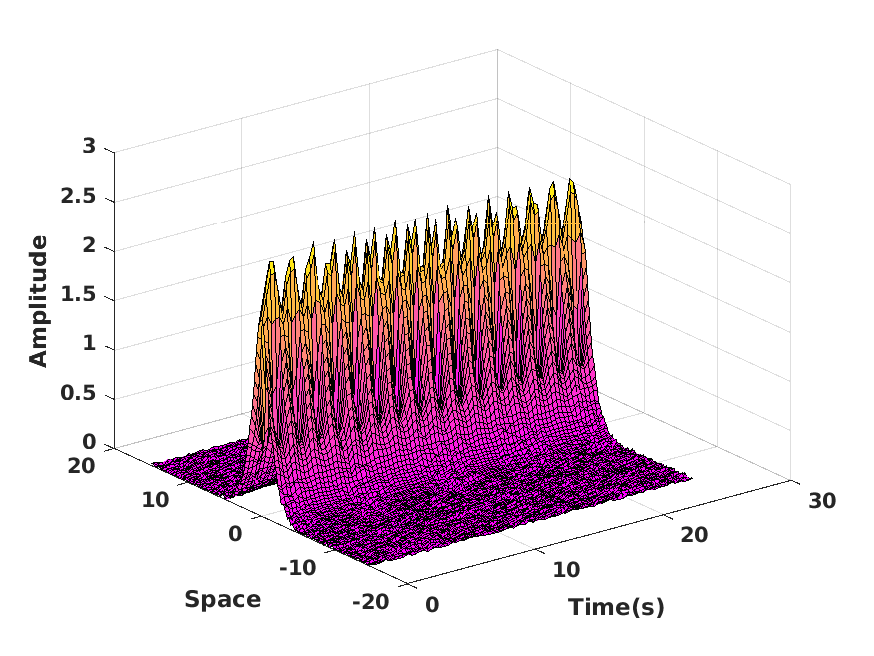}
    \end{adjustbox}
    \caption{ADM}
    \label{nlse_adm}
    \end{subfigure}
    \hspace{0.04\textwidth}
    \begin{subfigure}{0.47\textwidth}
    \begin{adjustbox}{max width=1\textwidth,center}
    \includegraphics{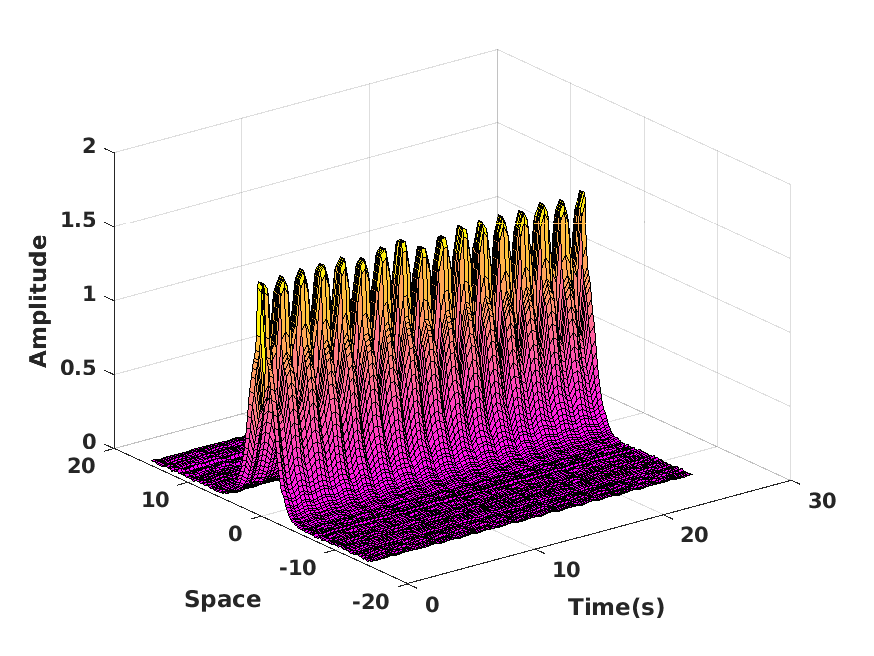}
    \end{adjustbox}
    \caption{Inexact ALM}
    \label{nlse_alm}
    \end{subfigure}
    \begin{subfigure}{0.47\textwidth}
    \begin{adjustbox}{max width=1\textwidth,center}
    \includegraphics{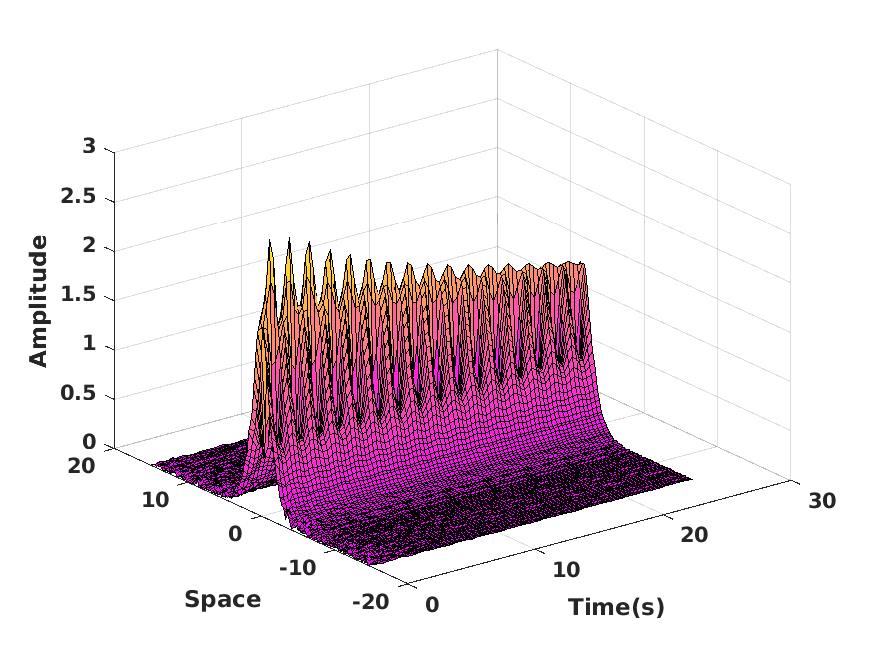}
    \end{adjustbox}
    \caption{TLS}
    \label{nlse_tls}
    \end{subfigure}
    \caption{NLSE}
    \label{nlse adm surf}
\end{figure}

\begin{figure}[H]
    \centering
    \begin{subfigure}{0.47\textwidth}
    \begin{adjustbox}{max width=1\textwidth,center}
    \includegraphics{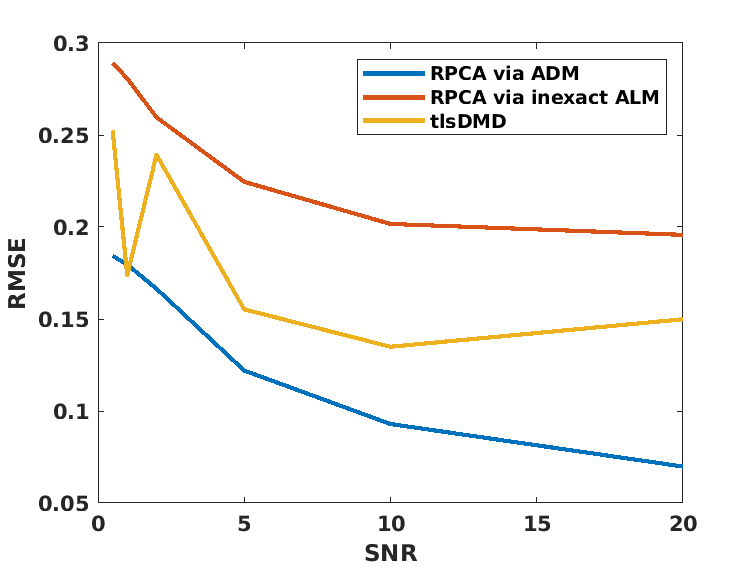}
    \end{adjustbox}
    \caption{RMSE}
    \label{RMSE_snr_nlse}
    \end{subfigure}
    \hspace{0.04\textwidth}
    \begin{subfigure}{0.47\textwidth}
    \begin{adjustbox}{max width=1\textwidth,center}
    \includegraphics{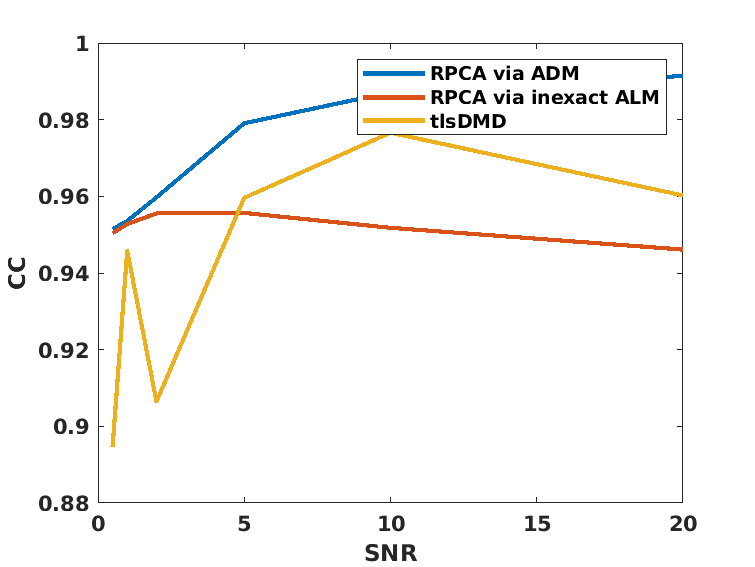}
    \end{adjustbox}
    \caption{CC}
    \label{CC_snr_nlse}
    \end{subfigure}
    \begin{subfigure}{0.47\textwidth}
    \begin{adjustbox}{max width=1\textwidth,center}
    \includegraphics{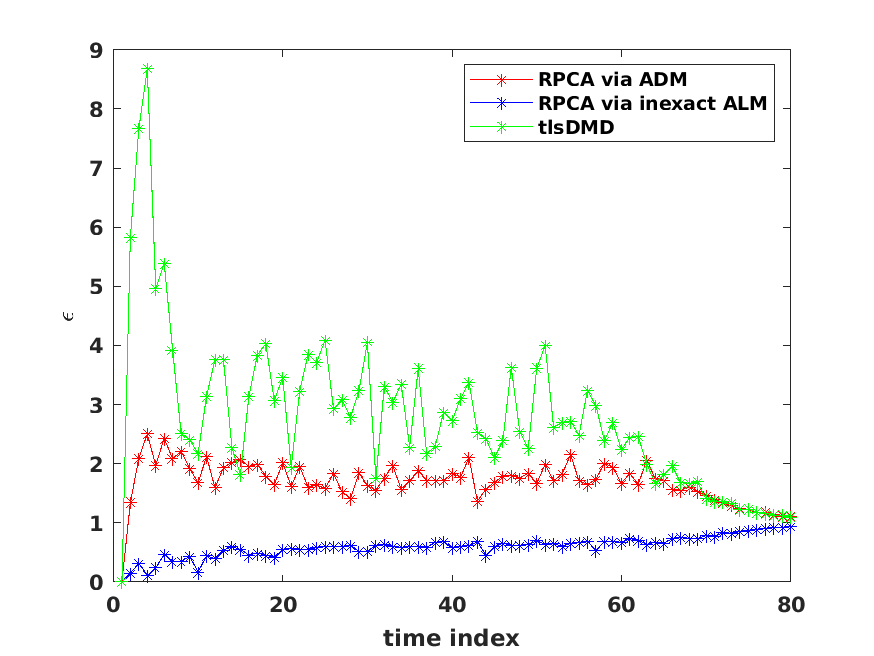}
    \end{adjustbox}
    \caption{Error variation with time}
    \end{subfigure}
    \begin{subfigure}{0.47\textwidth}
    \begin{adjustbox}{max width=1\textwidth,center}
    \includegraphics{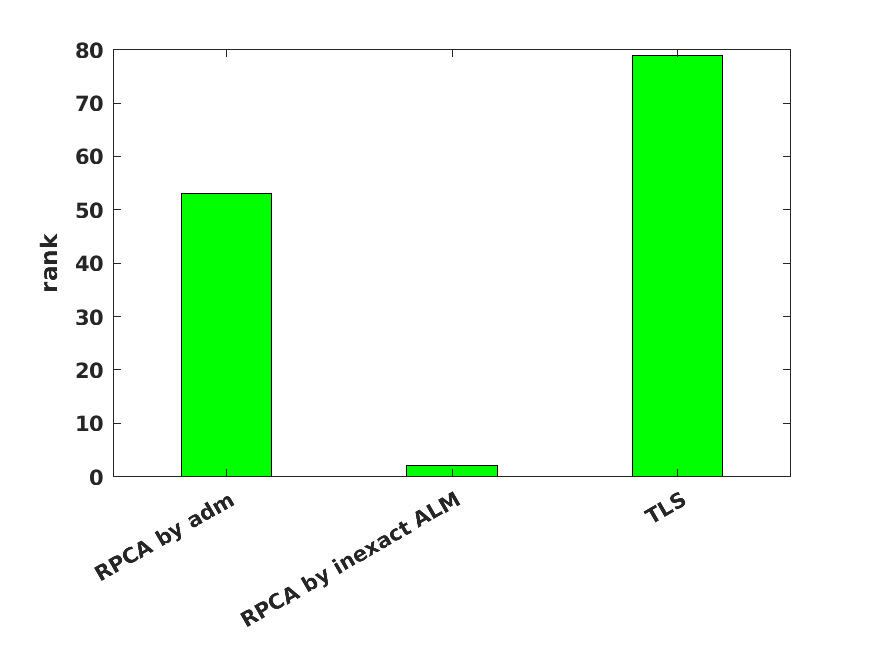}
    \end{adjustbox}
    \caption{Rank for each data-filtering method}
    \label{rank_nlse}
    \end{subfigure}
    \caption{NLSE}
    \label{nlse_error}
\end{figure}

The data for this example, Figure \ref{nlse_odat}, comes from the solution to the non-linear Schrodinger's equation (NLSE),
$$ \frac{\partial p}{\partial t} = \frac{i}{2} \frac{\partial^{2} p}{\partial w^2}  + i |p|^{2} p ,$$
where \(w\) and \(t\) are space and time, respectively. The domain are defined as:
$$ -15 < w < 15, $$
$$ 0 < t < 8 \pi. $$
The corrupted data (SNR = 20) is shown in Figure \ref{nlse_ndat}. The effect of RPCA via ADM, RPCA via inexact ALM and TLS-DMD on the noisy data is investigated and displayed in Figure \ref{nlse_adm}, Figure \ref{nlse_alm}, and Figure \ref{nlse_tls}, respectively. All of them tend to deal well with the noise in the data. The impact of the amount of the signal compared to that of noise is studied for each method based on RMSE and CC. It appears that RPCA via ADM performs much better that the other two methods for some given level of SNR, Figure \ref{RMSE_snr_nlse} and Figure \ref{CC_snr_nlse}. \\ \\
The approximation error in DMD is plotted in Figure \ref{nlse_error}. RPCA via inexact ALM sees a gradual increase in the error over time, whereas error in RPCA via ADM grows early on and then gets roughly constant. In case of the TLS-DMD, there is a sharp increase in the error in the very beginning after which a drop is witnessed.    

\begin{figure}[H]
    \centering
    \begin{subfigure}{0.47\textwidth}
    \begin{adjustbox}{max width=1\textwidth,center}
    \includegraphics{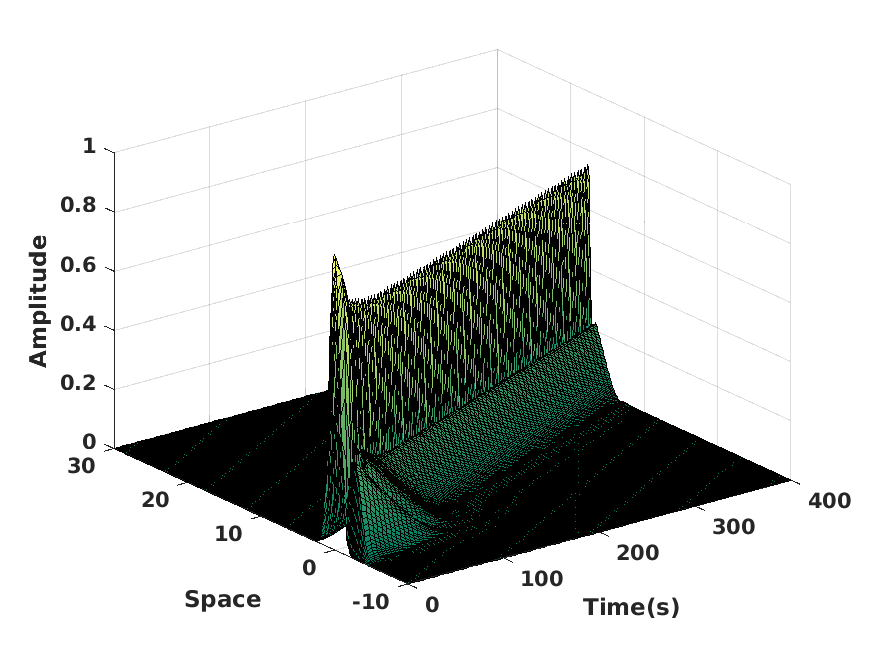}
    \end{adjustbox}
    \caption{Original Data}
    \label{FNE_odat}
    \end{subfigure}
    \hspace{0.04\textwidth}
    \begin{subfigure}{0.47\textwidth}
    \begin{adjustbox}{max width=1\textwidth,center}
    \includegraphics{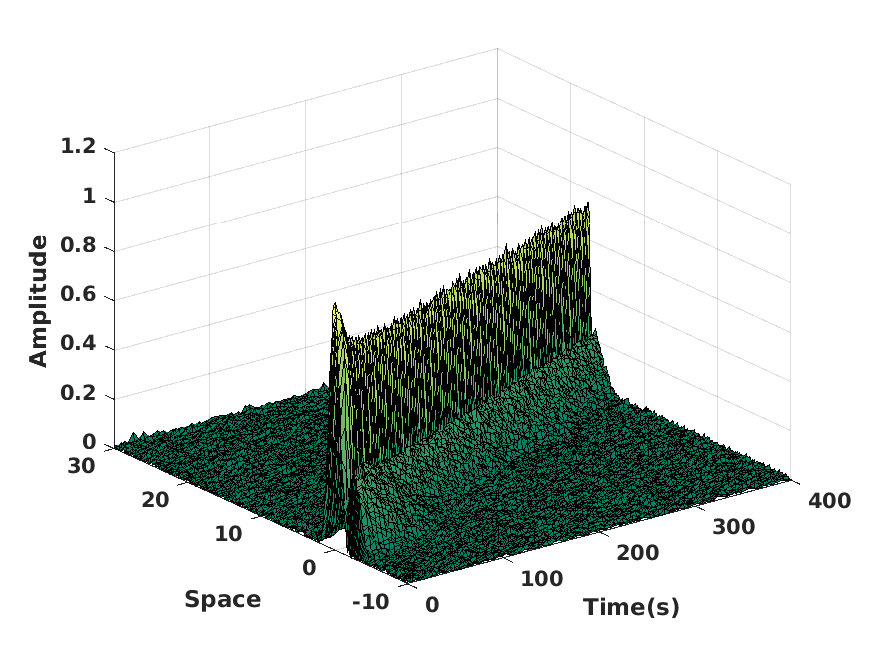}
    \end{adjustbox}
    \caption{Noisy Data}
    \label{FNE_ndat}
    \end{subfigure}
    \begin{subfigure}{0.47\textwidth}
    \begin{adjustbox}{max width=1\textwidth,center}
    \includegraphics{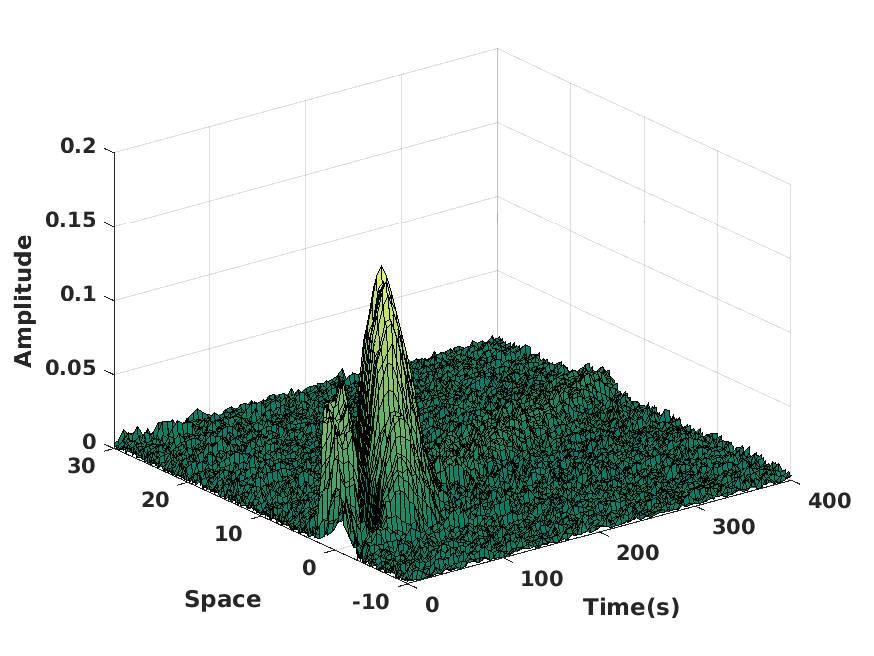}
    \end{adjustbox}
    \caption{ADM}
    \label{FNE_adm}
    \end{subfigure}
    \hspace{0.04\textwidth}
    \begin{subfigure}{0.47\textwidth}
    \begin{adjustbox}{max width=1\textwidth,center}
    \includegraphics{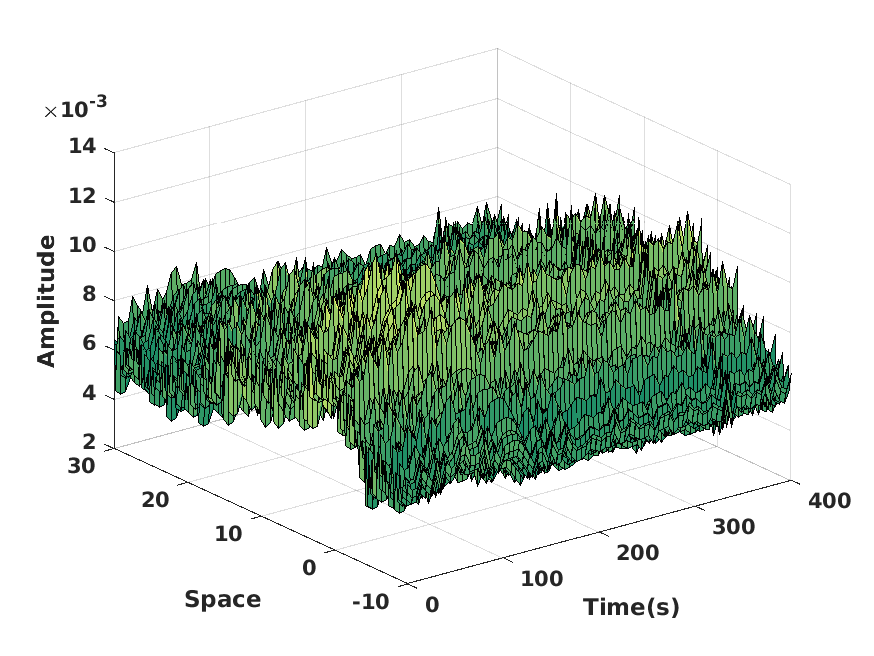}
    \end{adjustbox}
    \caption{Inexact ALM}
    \label{FNE_alm}
    \end{subfigure}
    \begin{subfigure}{0.47\textwidth}
    \begin{adjustbox}{max width=1\textwidth,center}
    \includegraphics{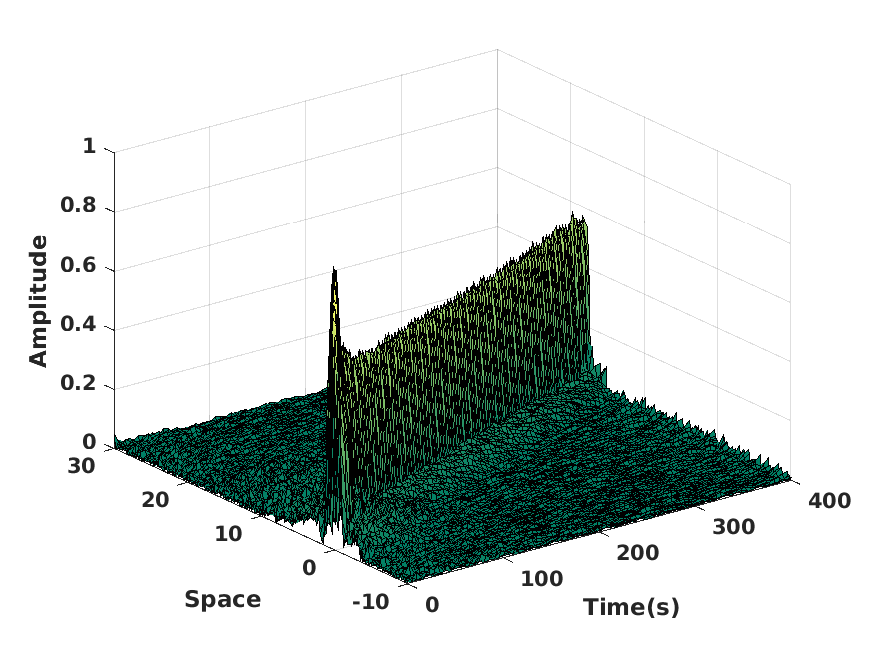}
    \end{adjustbox}
    \caption{TLS}
    \label{FNE_tls}
    \end{subfigure}
    \caption{FNE}
    \label{nlse adm surf}
\end{figure}

\begin{figure}[H]
    \centering
    \begin{subfigure}{0.47\textwidth}
    \begin{adjustbox}{max width=1\textwidth,center}
    \includegraphics{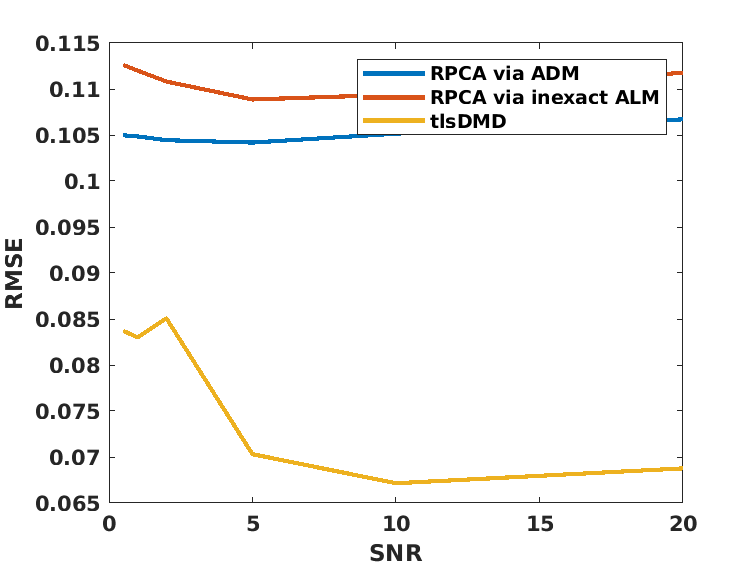}
    \end{adjustbox}
    \caption{RMSE}
    \label{fne_rmse}
    \end{subfigure}
    \hspace{0.04\textwidth}
    \begin{subfigure}{0.47\textwidth}
    \begin{adjustbox}{max width=1\textwidth,center}
    \includegraphics{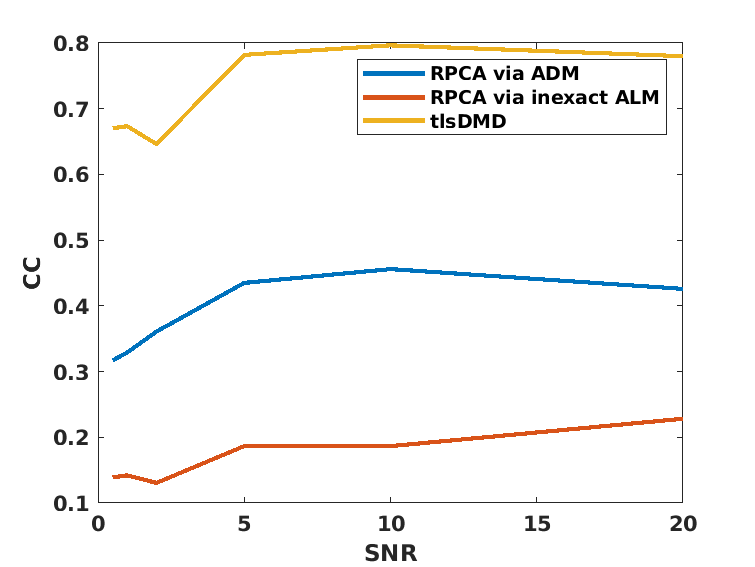}
    \end{adjustbox}
    \caption{CC}
    \label{fne_cc}
    \end{subfigure}
    \begin{subfigure}{0.47\textwidth}
    \begin{adjustbox}{max width=1\textwidth,center}
    \includegraphics{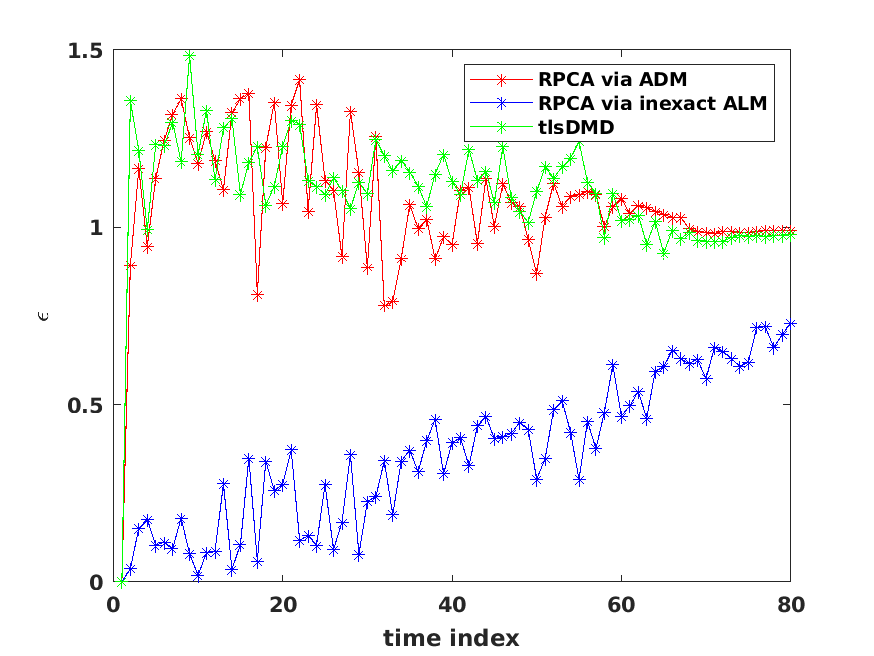}
    \end{adjustbox}
    \caption{Error variation with time}
    \label{error_fne}
    \end{subfigure}
    \begin{subfigure}{0.47\textwidth}
    \begin{adjustbox}{max width=1\textwidth,center}
    \includegraphics{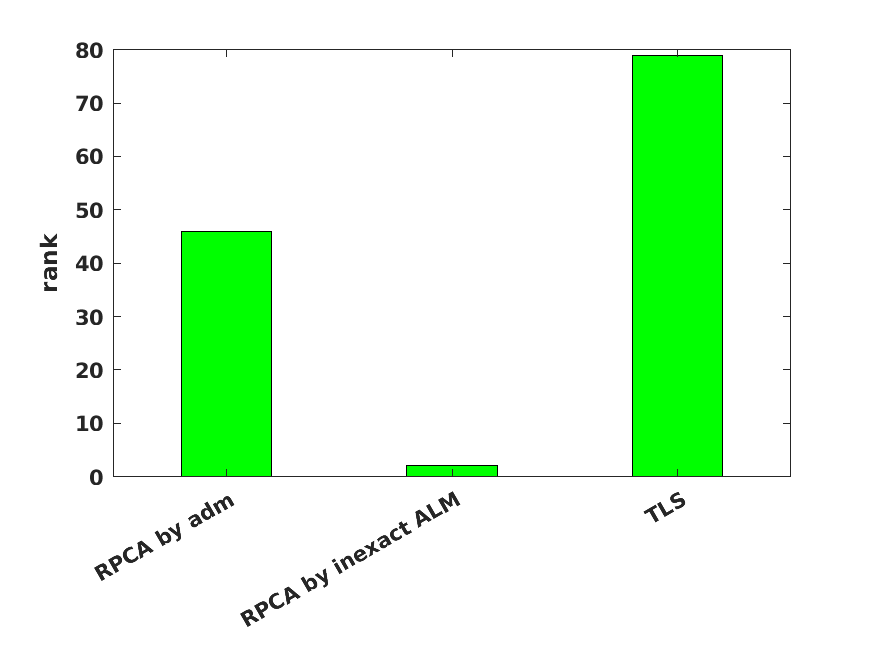}
    \end{adjustbox}
    \caption{Rank for each data-filtering method}
    \label{rank_fne}
    \end{subfigure}
    \caption{FNE}
    \label{nlse adm surf}
\end{figure}

\noindent \textbf{Fitzhugh-Nagumo Equation}\\ \\
The Fitzhugh-Nagumo equation (FNE) models the voltage spikes in neurons. FNE consists of the following system of partial differential equations,
\begin{equation}\label{FNE 1}
    \frac{\delta V}{\delta t} = D\frac{\delta^2 V}{\delta x^2} + V(a-V)(V-1) - W
\end{equation}
\begin{equation}\label{FNE 2}
    \frac{\delta W}{\delta t} = bV - cW
\end{equation}
where $V(x,t)$ measures the voltage in neuron. The initial conditions are as,
\begin{equation}\label{FNE 3}
    V(x,0) = exp(-x^2),
\end{equation}
\begin{equation}\label{FNE 4}
    W(x,0) = 0.2exp(-(x+2)^2).
\end{equation}
The boundary condition used on the partial differential equation system is,
\begin{equation} \label{FNE 5}
    \frac{\delta V}{\delta x} = 0 \;\; and \;\; \frac{\delta W}{\delta x} = 0 \;\; at\; x = a,b.
\end{equation}
The raw data and the noisy data (SNR of 20) are shown in Figure \ref{FNE_odat} and Figure \ref{FNE_ndat}, respectively. The noise correction is somewhat done by RPCA via ADM (Figure \ref{FNE_adm}), but it is poorly done by RPCA via inexact ALM Figure \ref{FNE_alm}. Among the three techniques, TLS-DMD happens to be the best in denoising the data-set, Figure \ref{FNE_tls}.\\ \\
RMSE is found to go down as SNR is increased for each method. The effect of the SNR is more pronounced on TLS-DMD than on RPCA via ADM or RPCA via inexact ALM, Figure \ref{fne_rmse} and Figure \ref{fne_cc}. The error in the DMD approximation goes up in the beginning and then stabilizes for both TLS-DMD and RPCA via ADM, whereas it is a slowly increasing trend for RPCA via inexact ALM, Figure \ref{error_fne}. 

\subsection{Shallow water equations}
\begin{figure}[H]
    \centering
    \begin{subfigure}{0.47\textwidth}
    \begin{adjustbox}{max width=1\textwidth,center}
    \includegraphics{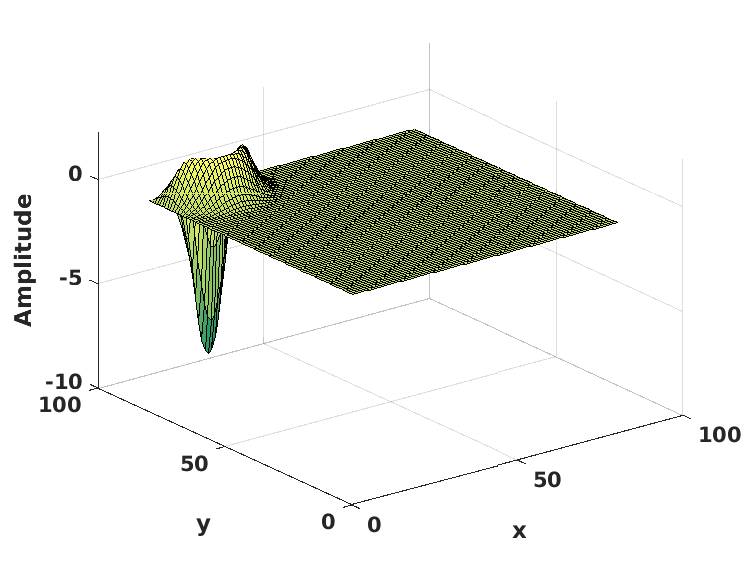}
    \end{adjustbox}
    \caption{Original Data}
    \label{sw_odat}
    \end{subfigure}
    \hspace{0.04\textwidth}
    \begin{subfigure}{0.47\textwidth}
    \begin{adjustbox}{max width=1\textwidth,center}
    \includegraphics{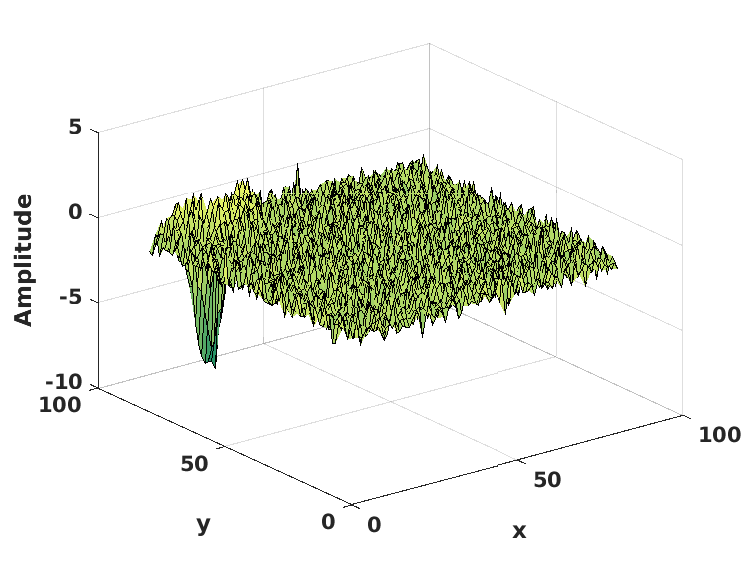}
    \end{adjustbox}
    \caption{Noisy Data}
    \label{sw_ndat}
    \end{subfigure}
    \begin{subfigure}{0.47\textwidth}
    \begin{adjustbox}{max width=1\textwidth,center}
    \includegraphics{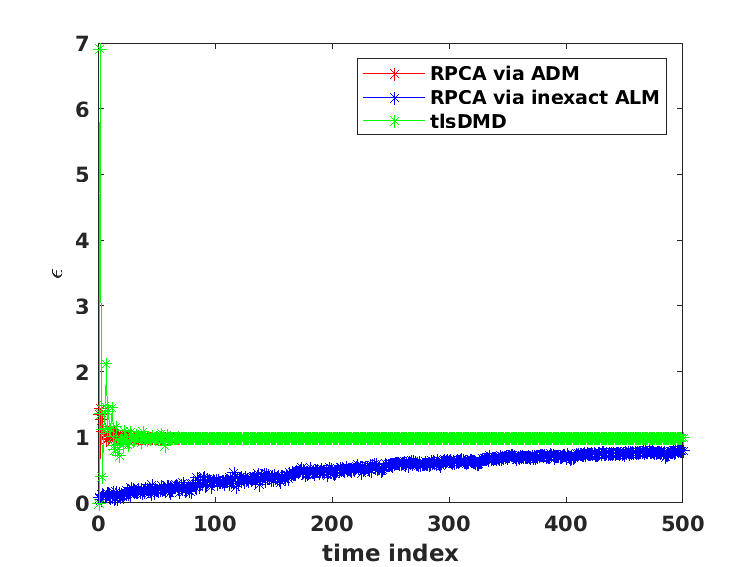}
    \end{adjustbox}
    \caption{Error variation with time}
    \label{sw_error}
    \end{subfigure}
    \begin{subfigure}{0.47\textwidth}
    \begin{adjustbox}{max width=1\textwidth,center}
    \includegraphics{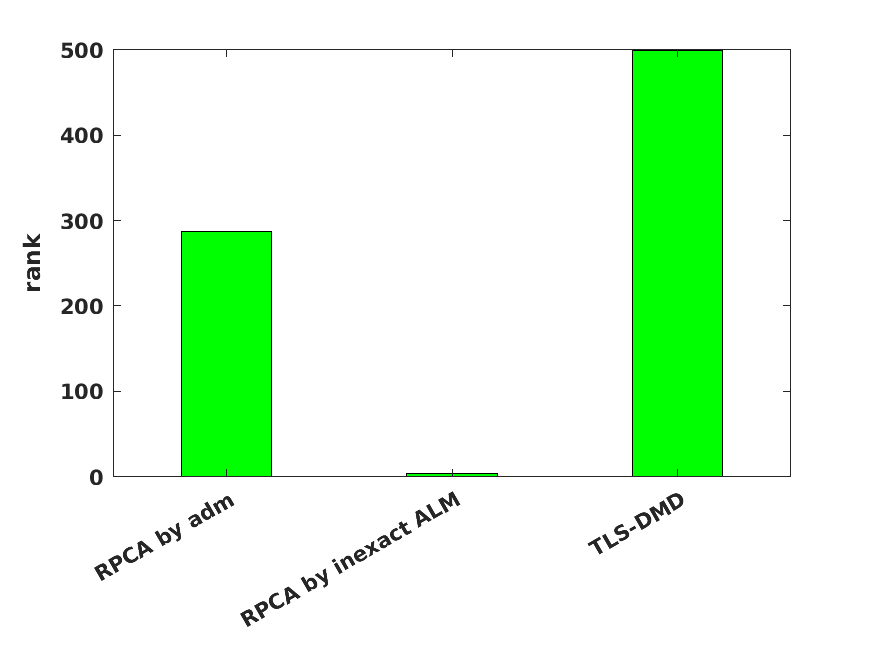}
    \end{adjustbox}
    \caption{Rank for each data-filtering method}
    \label{rank_sw}
    \end{subfigure}
    \hspace{0.04\textwidth}
    \caption{SWE}
\end{figure}

The motion of water in rivers and channels are often referred to as the shallow water equations (SWE) given by,
\begin{equation}
    \frac{\partial(\rho \kappa)}{\partial t} + \frac{\partial(\rho \kappa u)}{\partial x} + \frac{\partial(\rho \kappa v)}{\partial y} = 0,
\end{equation}

\begin{equation}
    \frac{\partial(\rho \kappa u)}{\partial t} + \frac{\partial}{\partial x} (\rho \kappa u^2 + \frac{1}{2} \rho g \kappa^2) + \frac{\partial (\rho \kappa u v)}{\partial y} = 0,
\end{equation}

\begin{equation}
    \frac{\partial(\rho \kappa v)}{\partial t} + \frac{\partial (\rho \kappa u v)}{\partial x} + \frac{\partial}{\partial y} (\rho \kappa v^2 + \frac{1}{2} \rho g \kappa^2) = 0,
\end{equation}
to model the dynamics of fluid column height (\(\kappa\)) of a fluid of constant density \(\rho\). Here, \((x,y)\) defines the horizontal 2D space, and \(u\) and \(v\) denote the flow velocity in the 2D space. The solution to SWE via finite differencing produces the data for this example. The raw data and the noisy data are illustrated in Figure \ref{sw_odat} and Figure \ref{sw_ndat}, respectively for a certain time instant (the tenth temporal node in this case). Much like the previous example, the uncertainty in the DMD approximation is similar for both RPCA via ADM and TLS-DMD, whereas it moderately increases with time for RPCA via ALM, Figure \ref{sw_error}. 

\section{Conclusion and Future Work}
This work is concerned with three major noise-correction tools that can be integrated with DMD. These tools were previously used in other applications like detecting fake moustache or glasses in an image. In this work, we exploit them to remove noise from the dataset generated from PDEs and eventually process the data via DMD to come up with reduced order models. The rank of the filtered data from each method is plotted in Figure \ref{rank_nlse}, \ref{rank_fne}, \ref{rank_sw}. It is high for both RPCA via ADM and TLS-DMD, but remains very low for RPCA via inexact ALM in all three examples. The effect of SNR on each method is examined and no absolute pattern is found. But, there exists a trend in the way the error in the DMD model accumulates over time. The uncertainty in the DMD model increases very early in time and then stabilizes to a certain value with time for RPCA via ADM and RPCA via inexact ALM method, whereas the error for TLS-DMD sees a gradual increase in time. This observations are made based on the results from testing the methods on three PDEs: the Non-linear Schrodinger Equation, the Fitzhugh-Nagumo Equation, and the Shallow water equations. In future, we would like to design a denoising scheme that can outperform the ones that are used in this work.
\label{CFW}

\section{Acknowledgement}
This project is partially funded by the Office of Research, North South University, under grant \textbf{CTRG-20/SEPS/08}.

\bibliographystyle{ieeetr}
\bibliography{dmd}
\end{document}